\documentclass[11pt, a4paper, twoside,leqno]{amsart}
\usepackage[centering, totalwidth = 420pt, totalheight = 610pt]{geometry}
\usepackage{amssymb, amsmath, amsthm, enumerate, microtype, stmaryrd, url,mathrsfs,mathpartir}
\usepackage[latin1]{inputenc}
\usepackage[dvips, arrow, matrix, tips, curve]{xy}
\usepackage{color}
\definecolor{darkgreen}{rgb}{0,0.45,0}
\usepackage[pagebackref,colorlinks,citecolor=darkgreen,linkcolor=darkgreen]{hyperref}
\SelectTips{cm}{10}

\newcommand{\cat}[1]{\mathbf{#1}}

\newcommand{\thg}{{\mathord{\text{--}}}}

\renewcommand{\c}{,\,\,}
\renewcommand{\t}{\,\mathbin{\vdash}\,}

\newcommand{\cdl}[2][]{\xymatrix@1#1{#2}}

\newcommand{\st}{\ \mathsf{sort}}
\newcommand{\ty}{\ \mathsf{type}}
\newcommand{\Id}{\mathsf{Id}}

\newcommand{\J}{{\mathcal J}}

\newcommand{\xtor}[1]{\cdl[@1]{{} \ar[r]|-{\object@{|}}^{#1} & {}}}

\makeatletter
\def\mathrlap@#1#2{\rlap{$\m@th#1{#2}$}}
\providecommand*\mathrlap[1][\@empty]{
  \ifx\@empty#1\@empty
    \expandafter \mathpalette \expandafter \mathrlap@
  \else
    \expandafter \MT_mathrlap:Nn \expandafter #1
  \fi
}
\def\hookleftarrowfill@{\arrowfill@\leftarrow\relbar{\relbar\joinrel\rhook}}
\def\twoheadleftarrowfill@{\arrowfill@\twoheadleftarrow\relbar\relbar}
\def\leftbararrowfill@{\arrowdoublefill@{\leftarrow\mkern-5mu}\relbar\mapstochar\relbar\relbar}
\def\Leftbararrowfill@{\arrowdoublefill@{\Leftarrow\mkern-2mu}\Relbar\Mapstochar\Relbar\Relbar}
\def\leftringarrowfill@{\arrowdoublefill@{\leftarrow\mkern-3mu}\relbar{\mkern-3mu\circ\mkern-2mu}\relbar\relbar}
\def\lefttriarrowfill@{\arrowfill@{\mathrel\triangleleft\mkern0.5mu\joinrel\relbar}\relbar\relbar}
\def\Lefttriarrowfill@{\arrowfill@{\mathrel\triangleleft\mkern1mu\joinrel\Relbar}\Relbar\Relbar}

\def\hookrightarrowfill@{\arrowfill@{\lhook\joinrel\relbar}\relbar\rightarrow}
\def\twoheadrightarrowfill@{\arrowfill@\relbar\relbar\twoheadrightarrow}
\def\rightbararrowfill@{\arrowdoublefill@{\relbar\mkern-0.5mu}\relbar\mapstochar\relbar\rightarrow}
\def\Rightbararrowfill@{\arrowdoublefill@{\Relbar\mkern-2mu}\Relbar\Mapstochar\Relbar\Rightarrow}
\def\rightringarrowfill@{\arrowdoublefill@\relbar\relbar{\mkern-2mu\circ\mkern-3mu}\relbar{\mkern-3mu\rightarrow}}
\def\righttriarrowfill@{\arrowfill@\relbar\relbar{\relbar\joinrel\mkern0.5mu\mathrel\triangleright}}
\def\Righttriarrowfill@{\arrowfill@\Relbar\Relbar{\Relbar\joinrel\mkern1mu\mathrel\triangleright}}

\def\leftrightarrowfill@{\arrowfill@\leftarrow\relbar\rightarrow}
\def\mapstofill@{\arrowfill@{\mapstochar\relbar}\relbar\rightarrow}

\renewcommand*\xleftarrow[2][]{\ext@arrow 20{20}0\leftarrowfill@{#1}{#2}}
\providecommand*\xLeftarrow[2][]{\ext@arrow 60{22}0{\Leftarrowfill@}{#1}{#2}}
\providecommand*\xhookleftarrow[2][]{\ext@arrow 10{20}0\hookleftarrowfill@{#1}{#2}}
\providecommand*\xtwoheadleftarrow[2][]{\ext@arrow 60{20}0\twoheadleftarrowfill@{#1}{#2}}
\providecommand*\xleftbararrow[2][]{\ext@arrow 10{22}0\leftbararrowfill@{#1}{#2}}
\providecommand*\xLeftbararrow[2][]{\ext@arrow 50{24}0\Leftbararrowfill@{#1}{#2}}
\providecommand*\xleftringarrow[2][]{\ext@arrow 10{26}0\leftringarrowfill@{#1}{#2}}
\providecommand*\xlefttriarrow[2][]{\ext@arrow 80{24}0\lefttriarrowfill@{#1}{#2}}
\providecommand*\xLefttriarrow[2][]{\ext@arrow 80{24}0\Lefttriarrowfill@{#1}{#2}}

\renewcommand*\xrightarrow[2][]{\ext@arrow 01{20}0\rightarrowfill@{#1}{#2}}
\providecommand*\xRightarrow[2][]{\ext@arrow 04{22}0{\Rightarrowfill@}{#1}{#2}}
\providecommand*\xhookrightarrow[2][]{\ext@arrow 00{20}0\hookrightarrowfill@{#1}{#2}}
\providecommand*\xtwoheadrightarrow[2][]{\ext@arrow 03{20}0\twoheadrightarrowfill@{#1}{#2}}
\providecommand*\xrightbararrow[2][]{\ext@arrow 01{22}0\rightbararrowfill@{#1}{#2}}
\providecommand*\xRightbararrow[2][]{\ext@arrow 04{24}0\Rightbararrowfill@{#1}{#2}}
\providecommand*\xrightringarrow[2][]{\ext@arrow 01{26}0\rightringarrowfill@{#1}{#2}}
\providecommand*\xrighttriarrow[2][]{\ext@arrow 07{24}0\righttriarrowfill@{#1}{#2}}
\providecommand*\xRighttriarrow[2][]{\ext@arrow 07{24}0\Righttriarrowfill@{#1}{#2}}

\providecommand*\xmapsto[2][]{\ext@arrow 01{20}0\mapstofill@{#1}{#2}}
\providecommand*\xleftrightarrow[2][]{\ext@arrow 10{22}0\leftrightarrowfill@{#1}{#2}}
\providecommand*\xLeftrightarrow[2][]{\ext@arrow 10{27}0{\Leftrightarrowfill@}{#1}{#2}}

\makeatother

\newcommand{\twocong}[2][0.5]{\ar@{}[#2] \save ?(#1)*{\cong}\restore}
\newcommand{\twoeq}[2][0.5]{\ar@{}[#2] \save ?(#1)*{=}\restore}
\newcommand{\rtwocell}[3][0.5]{\ar@{}[#2] \ar@{=>}?(#1)+/l 0.2cm/;?(#1)+/r 0.2cm/^{#3}}
\newcommand{\ltwocell}[3][0.5]{\ar@{}[#2] \ar@{=>}?(#1)+/r 0.2cm/;?(#1)+/l 0.2cm/^{#3}}
\newcommand{\ltwocello}[3][0.5]{\ar@{}[#2] \ar@{=>}?(#1)+/r 0.2cm/;?(#1)+/l 0.2cm/_{#3}}
\newcommand{\dtwocell}[3][0.5]{\ar@{}[#2] \ar@{=>}?(#1)+/u  0.2cm/;?(#1)+/d 0.2cm/^{#3}}
\newcommand{\dthreecell}[3][0.5]{\ar@{}[#2] \ar@3{->}?(#1)+/u  0.2cm/;?(#1)+/d 0.2cm/^{#3}}
\newcommand{\utwocell}[3][0.5]{\ar@{}[#2] \ar@{=>}?(#1)+/d 0.2cm/;?(#1)+/u 0.2cm/_{#3}}
\newcommand{\dtwocelltarg}[3][0.5]{\ar@{}#2 \ar@{=>}?(#1)+/u  0.2cm/;?(#1)+/d 0.2cm/^{#3}}
\newcommand{\utwocelltarg}[3][0.5]{\ar@{}#2 \ar@{=>}?(#1)+/d  0.2cm/;?(#1)+/u 0.2cm/_{#3}}

\newdir{ >}{{}*!/-5pt/\dir{>}}

\makeatletter
\newcommand{\pgph}{\paragraph{\!}}
\renewcommand{\paragraph}{\@startsection
{paragraph}%
{3}%
{0mm}%
{-\baselineskip}%
{-0.4em plus 0.2em minus 0.2em}%
{\normalfont\normalsize\bfseries}}%
\newcommand{\paragraphprime}{\@startsection
{paragraph}%
{4}%
{0mm}%
{-\baselineskip}%
{-0.4em plus 0.2em minus 0.2em}%
{\normalfont\normalsize\bfseries}}%
\makeatother

\numberwithin{equation}{section} \numberwithin{paragraph}{section}

\newenvironment{Thm}{\paragraph{Theorem:}\em}{\vskip\baselineskip}

\newenvironment{Prop}{\paragraph{Proposition:}\em}{\vskip\baselineskip}

\newenvironment{Exs}{\paragraph{Examples:}}{\vskip\baselineskip}

\newenvironment{Lems}{\subsubsection*{\textbf{\upshape{Lemma}}}\em}{\vskip\baselineskip}

\makeatletter \@namedef{itemize*}{\itemize\parsep\z@ \parskip\z@}
\@namedef{enditemize*}{\enditemize} \@namedef{enumerate*}{\enumerate\parsep\z@
\parskip\z@} \@namedef{endenumerate*}{\endenumerate}
\def\Pr@@f{\subsubsection*{\textbf{Proof}}}
\def\pr@@f[#1]{\subsubsection*{{\textbf{Proof}} #1}}
\makeatother

\makeatletter \@namedef{itemize*}{\itemize\parsep\z@ \parskip\z@}
\@namedef{enditemize*}{\enditemize} \@namedef{enumerate*}{\enumerate\parsep\z@
\parskip\z@} \@namedef{endenumerate*}{\endenumerate} \makeatother

\makeatletter
\def\matrixobject@{%
  \edef \next@{={\DirectionfromtheDirection@ }}%
  \expandafter \toks@ \next@ \plainxy@
  \let\xy@@ix@=\xyq@@toksix@
  \xyFN@ \OBJECT@}
\let\xy@entry@@norm=\entry@@norm
\def\entry@@norm@patched{%
  \let\object@=\matrixobject@
  \xy@entry@@norm }
\AtBeginDocument{\let\entry@@norm\entry@@norm@patched}
\makeatother

\newcommand{\defeq}{\mathrel{\mathop:}=}

\begin{document}
\title[On the strength of dependent products]{On the strength of dependent products\\in the type theory of Martin-Löf}
\author{Richard Garner}
\address{Department of Mathematics, Uppsala University, Box 480, S-751 06 Uppsala, Sweden}
\email{rhgg2@cam.ac.uk}
\thanks{Supported by a Research Fellowship of St John's College, Cambridge and a Marie Curie Intra-European Fellowship, Project No.\ 040802.}
\begin{abstract}
One may formulate the dependent product types of Martin-L\"of type theory
either in terms of abstraction and application operators like those for the
lambda-calculus; or in terms of introduction and elimination rules like those
for the other constructors of type theory. It is known that the latter rules
are at least as strong as the former: we show that they are in fact strictly
stronger. We also show, in the presence of the identity types, that the
elimination rule for dependent products~--~which is a ``higher-order''
inference rule in the sense of Schroeder-Heister~--~can be reformulated in a
first-order manner. Finally, we consider the principle of function
extensionality in type theory, which asserts that two elements of a dependent
product type which are pointwise propositionally equal, are themselves
propositionally equal. We demonstrate that the usual formulation of this
principle fails to verify a number of very natural propositional equalities;
and suggest an alternative formulation which rectifies this deficiency.
\end{abstract}

 \maketitle 
\section{Introduction} This is the first in a series of papers recording the
author's investigations into the semantics of Martin-L\"of's dependent type
theory; by which we mean the type theory set out in the expository
volume~\cite{Programming}. The main body of these investigations concerns what
the author is calling \emph{two-dimensional} models of dependent type theory.
Recall that one typically divides the models of Martin-L\"of's type theory into
\emph{extensional} and \emph{intensional} ones, the former differentiating
themselves from the latter by their admission of an equality reflection rule
which collapses the propositional and definitional equalities of the language
into a single, \emph{judgemental}, equality. The two-dimensional models that
the author is studying are of the intensional kind, but are not wholly
intensional: they admit instances of the equality reflection rule at just those
types which are themselves identity types.

In the process of making his investigations, the author has discovered certain
unresolved issues concerning the dependent product types of Martin-L\"of type
theory; and since these issues exist beyond the domain of two-dimensional
models, it seemed worthwhile to collect his conclusions into this preliminary
paper.

The first of these issues concerns how we formulate of the rules for the
dependent product types. There are two accepted ways of doing this. In both
cases, we begin with a \emph{formation} rule which, given a type $A$ and a type
$B(x)$ dependent on $x : A$, asserts the existence of a type $\Pi(A, B)$; and
an \emph{abstraction} rule which says that, from an element $f(x) : B(x)$
dependent on $x : A$, we may deduce the existence of an element
\mbox{$\lambda(f) : \Pi(A, B)$}. We may then complement these rules with either
an \emph{application} rule, which tells us that, from $m : \Pi(A, B)$ and $a :
A$, we may infer an element \mbox{$\mathsf{app}(m, x) : B(x)$}; or an
\emph{elimination} rule, which essentially tells us that any (dependent) function out
of $\Pi(A, B)$ is determined, up-to-propositional-equality, by its values on
those elements of the form $\lambda(f)$ for some dependent element $x : A \t
f(x) : B(x)$.

There are two problematic features here. The first concerns the nature of the
elimination rule, which is a \emph{higher-order inference rule} in the sense of
Schroeder-Heister~\cite{PSH}. In order to formulate this rule rigorously, we must
situate our type theory within an ambient calculus possessing higher-order
features; a suitable choice being the \emph{Logical Framework} described in
Part~III of~\cite{Programming}, and recalled in Section~\ref{typeinlf} below.
Yet it may be that we do not wish to do this: one reason being that the
categorical semantics of Martin-L\"of type theory looks rather different when
it is formulated within the Logical Framework. Hence our first task in this
paper is to give a first-order reformulation of the elimination rule in terms
of the application rule and a propositional form of the $\eta$-rule; a
reformulation that may be stated without recourse to the Logical Framework.

The second problematic feature concerns the precise relationship between the
application and elimination rules for dependent products. We know that the
application rule may be defined in terms of the elimination rule, so that the
elimination rule is stronger; yet it is not known whether it is \emph{strictly}
stronger. Our second task is to show that this is in fact the case; we do this
by describing a non-standard interpretation of the $\Pi$-types for which the
application rule obtains, yet not the elimination rule.

We then move on to another issue, namely the formulation of the principle of \emph{function
extensionality} in Martin-L\"of type theory. This principle asserts that if $m$ and $n$ are
elements of $\Pi(A, B)$ and we can affirm a propositional equality between $\mathsf{app}(m,
x)$ and $\mathsf{app}(n, x)$ whenever $x : A$, then we may deduce the existence of a
propositional equality between $m$ and $n$.  One result of the author's investigations has
been that, if we are to obtain a notion of two-dimensional model which is reasonably urbane
from a category-theoretic perspective, then we must impose some kind of function
extensionality. Yet the principle of function extensionality just stated has been found
wanting in this regard, since it fails to provide witnesses for a number of very natural
propositional equalities which are demanded by the  semantics; some of which are detailed in
Examples \ref{section56} below. From a category-theoretic perspective, we might say that the
principle of function extensionality fails to be \emph{coherent}. Our third task in this
paper, therefore, is to propose a suitably coherent replacement for function extensionality.

\textbf{Acknowledgements}. The author wishes to thank Johan Granstr\"om, Per
Martin-L\"of, Erik Palmgren, Olov Wilander and other members of the
Stockholm-Uppsala Logic Seminar for useful comments and suggestions on earlier
drafts of this paper. He also thanks the anonymous referee for several useful
comments and suggestions.

\section{Martin-L\"of type theory}\label{typeinlf}
\pgph We begin with a brief summary of the two principal ways in which one may
present Martin-L\"of type theory. The more straightforward is the
``polymorphic'' presentation of \cite{pml, pml2}. This is given by a reasoning system with four
basic forms of judgement:
\begin{itemize*}
 \item $\Gamma \vdash A \ty$ (``$A$ is a type under the hypothesis $\Gamma$'');
 \item $\Gamma \vdash a : A$ (``$a$ is an element of $A$ under the hypothesis $\Gamma$);
 \item $\Gamma \vdash A = B \ty$ (``$A$ and $B$ are equal types under the hypothesis $\Gamma$'');
 \item $\Gamma \vdash a = b : A$ (``$a$ and $b$ are equal elements of $A$ under the
 hypothesis
 $\Gamma$'').
\end{itemize*}
Here, $\Gamma$ is to be a \emph{context} of assumptions, $\Gamma = (x_1 : A_1\c
x_2 : A_2\c \dots\c x_n : A_{n-1})$, subject to a requirement of
well-formedness which affirms that each $A_i$ is a type under the assumptions
$(x_1 : A_1\c \dots\c x_{i-1} : A_{i-1})$. The polymorphic presentation of
Martin-L\"of type theory is now given by specifying a sequent calculus over
these four forms of judgement: so a number of \emph{axiom} judgements, together
with a number of \emph{inference rules}
\begin{equation*}
\inferrule{\J_1 \\ \cdots \\ \J_n}{\J}
\end{equation*}
allowing us to derive the validity of the judgement $\J$ from that of the
$\J_i$'s. As usual, these inference rules separate into a group of
\emph{structural rules} which deal with the contextual book-keeping of
weakening, contraction, exchange and substitution; and a group of \emph{logical
rules}, which describe the constructions we wish to be able to perform inside
our type theory: constructions such as cartesian product of types, disjoint
union of types, or formation of identity types.

\pgph However, the polymorphic presentation of type theory will be inadequate
for our purposes, because the elimination rule for dependent products we wish
to study requires the use of a \emph{second-order judgement} $\Gamma \vdash
\mathcal J$, in which the context of assumptions $\Gamma$ itself contains a
judgement under hypotheses. One solution to this problem is suggested by
Troelstra and van Dalen in~\cite[Chapter 11]{IIM}: we extend our system with
explicit second-order judgement forms expressing that ``$B$ is a family of
types over $A$ under the hypothesis $\Gamma$'', and so on, and express the
elimination rule in terms of these.
%
A second solution---and the one we adopt here---makes use of the
``monomorphic'' presentation of Martin-L\"of type theory. This is given in
terms of the \emph{Logical Framework}, which is essentially a formalisation of
the meta-theory we use to reason about the calculus of types. The basic
judgements of this meta-theory look rather like those of type theory:
\begin{mathpar}
 \Gamma \vdash A \st\text;

 \Gamma \vdash a : A\text;

 \Gamma \vdash A = B \st\text;

 \text{and} \quad \Gamma \vdash a = b : A\text.
\end{mathpar}
However, the meaning is somewhat different. We think of a sort of the Logical
Framework as being a category of \emph{judgements about type theory}. In
particular, the Logical Framework has rules
\begin{equation*}
    \inferrule { }
    {\t \textrm{type} \st} \qquad \text{and} \qquad
    \inferrule { }
    {A : \textrm{type} \t \textrm{el}\,A \st}\text,
\end{equation*}
which express the existence of the category of judgements ``--- is a type'';
and, under the assumption that ``$A$ is a type'', of the category of judgements
``--- is an element of $A$''. Using these, we may interpret more complex
judgements of type theory; for example, if we know that ``$A$ is a type'', then
we can interpret the judgement $\J$ that ``$B(x)$ is a type under the
hypothesis that $x$ is an element of $A$'' as
\begin{equation*}
    x : \textrm{el}\,A \t B(x) : \text{type}\text.
\end{equation*}
Yet this is not an entirely faithful rendition of $\J$, since strictly
speaking, the displayed sequent asserts the judgement ``$B(x)$ is a type''
under the hypothesis that ``$x$ is an element of $A$''. To resolve this, we
introduce the other key aspect of the Logical Framework, namely the
\emph{function sorts}. These are specified by rules of formation, abstraction
and application:
\begin{mathpar}
    \inferrule{\Gamma\c x : A \t B (x) \st}{\Gamma \t (x : A)\,B \st}\text,

    \inferrule{\Gamma\c x : A \t b (x)  : B (x) }{\Gamma \t [x : A]\, b  (x) : (x : A)\,B (x) }

    \text{and} \qquad \inferrule{\Gamma\t f : (x : A)\,B \\ \Gamma \t a : A}{\Gamma \t f(a) : B(a)}
\end{mathpar}
subject to the $\alpha$-, $\beta$-, $\eta$-\ and $\xi$-rules of the lambda-calculus. Using function sorts, we can now render the judgement $\J$ more
correctly. We have the sort $(x : \textrm{el}\,A)\,\text{type}$, which is the
category of judgements ``--- is a type under the hypothesis that $x$ is an
element of $A$''; and can now interpret $\J$ as the judgement
\begin{equation*}
    \t B : (x : \textrm{el}\, A)\,\text{type}\text.
\end{equation*}

\pgph \label{translate} We may translate the polymorphic presentation of
Martin-L\"of type theory into the monomorphic one by encoding the inference
rules of the former as higher-order judgements of the latter. For instance,
consider the hypothetical type constructor $\Phi$ with rules
\begin{equation*}
    \inferrule{A \ty}{\Phi(A) \ty}\qquad \text{and} \qquad
    \inferrule{A \ty \\ a : A}{\phi_A(a) : \Phi(A)}\text.
\end{equation*}
We may encode this in the Logical Framework by terms
\begin{equation*}
    \t \Phi : (A : \text{type})\, \text{type} \qquad \text{and} \qquad
    \t \phi : (A : \text{type}\c a : \text{el}\, A)\,\text{el}\, \Phi(A)\text,
\end{equation*}
where for readability we write iterated function spaces as $(A : \text{type}\c
a : \text{el}\, A)\,\text{el}\, \Phi(A)$ instead of the more correct $(A :
\text{type})(a : \text{el}\, A)\,\text{el}\, \Phi(A)$. Note that this encoding
says more than the original, by affirming a certain insensitivity to ambient
context; since from the constants $\Phi$ and $\phi$, we obtain a whole family
of inference rules
\begin{equation*}
    \inferrule{\Gamma \t A \ty}{\Gamma \t \Phi(A) \ty}\qquad \text{and} \qquad
    \inferrule{\Gamma \t A \ty \\ \Gamma \t a : A}{\Gamma \t \phi_A(a) :
    \Phi(A)}\text,
\end{equation*}
together with further rules expressing stability under substitution in
$\Gamma$. However, this is no bad thing, since any acceptable inference rule of
the polymorphic theory will possess this ``naturality'' in the context
$\Gamma$. In the remainder of this paper we work in the monomorphic
presentation of type theory, but will take advantage of the above encoding
process in order to present the rules of our type theory in the more readable
polymorphic style. For more on the relationship between the monomorphic and polymorphic presentations, see \cite{HTT}.

\section{A first-order reformulation of the $\Pi$-elimination rule}
\pgph Our main concern in this paper is with the dependent product types of
Martin-L\"of type theory: but in this analysis, we will from time to time make
use of the \emph{identity types}, which are a reflection at the type level of
the equality judgements $a = b : A$. We begin, therefore, by recalling the
rules for the identity types:

\begin{equation*}
\inferrule*[right=$\Id$-form;]{A \ty \\ a, b : A}{\Id_A(a, b) \ty} \qquad
\inferrule*[right=$\Id$-intro;]{A \ty \\
a : A}{r(a) : \Id_A(a, a)}
\end{equation*}
\begin{equation*}
\inferrule*[right=$\Id$-elim;]{
  A \ty \\
  x, y : A \c z : \Id_A(x, y) \t C(x, y, z) \ty \\
  x : A \t d(x) : C(x, x, r(x))\\
  a, b : A \\
  p : \Id_A(a, b)
  }{J(d, a, b, p) : C(a, b, p)}
\end{equation*}
\begin{equation*}
\inferrule*[right=$\Id$-\textsc{comp}.]{
  A \ty \\
  x, y : A \c z : \Id_A(x, y) \t C(x, y, z) \ty \\
  x : A \t d(x) : C(x, x, r(x))\\
  a : A
  }{J(d, a, a, r(a)) = d(a) : C(a, a, r(a))}
\end{equation*}

The notion of equality captured by the identity types is known as \emph{propositional
equality}: to say that $a$ and $b$ are propositionally equal as elements of $A$ is to say that
we may affirm a judgement $p : \Id_A(a, b)$. We think of $\Id_A$ as being a type inductively
generated by the elements $r(a) : \Id_A(a, a)$, with the elimination rule and computation
rules expressing that any dependent function out of $\Id_A$ is determined
up-to-propositional-equality by its value on elements of the form $r(a)$.
%
%

\pgph \label{appform} We are now ready to describe the two standard formulations of dependent
product types in Martin-L\"of type theory. The first, which we will refer to as the
$\mathsf{app}$-formulation, is analogous to the $\lambda$-calculus with the $\beta$-rule but
no $\eta$-rule:
\begin{mathpar}
\inferrule*[right=$\Pi$-form;]{A \ty \\ x : A \t B(x) \ty}{\Pi(A, B) \ty}

\inferrule*[right=$\Pi$-abs;]{
x : A \t f(x) : B(x)}{\lambda(f) : \Pi(A, B)}

\inferrule*[right=$\Pi$-app;]{
  m : \Pi(A, B) \\
  a : A}{
  \mathsf{app}(m, a) : B(a)}

\inferrule*[right=$\Pi$-$\beta$.]{
  x : A \t f(x) : B(x) \\
  a : A}{
  {\mathsf{app}(\lambda(f), a) = f(a) : B(a)}}
\end{mathpar}
Note that, for the sake of readability we omit the hypotheses $A \ty$ and $x : A \t
B(x) \ty$ from the last three of these rules; and in future, we may omit any
such hypotheses that are reconstructible from the context. To further reduce
syntactic clutter, we may also write $\Pi x : A.\, B(x)$ instead of $\Pi(A, [x
: A]\,B(x))$; $\lambda x.\, f(x)$ instead of $\lambda([x : A]\, f(x))$; and
$m \cdot a$ instead of $\mathsf{app}(m, a)$.


\pgph \label{funsplitform} As we noted in the Introduction, the second formulation of
dependent products---which we will refer to as the \textsf{funsplit}-formulation---has the
same introduction and abstraction rules but replaces the application and $\beta$-rules with
elimination and computation rules which mirror those for the other constructors of type
theory: they assert that each type $\Pi(A, B)$ is inductively generated by the elements of the
form $\lambda(f)$.
\begin{equation*}
\inferrule*[right=$\Pi$-elim;]{
  y : \Pi(A, B) \t C(y) \ty \\
  f : (x : A)\, B(x) \t d(f) : C(\lambda(f)) \\
  m : \Pi(A, B)}{
  \textsf{funsplit}(d, m) : C(m)}
\end{equation*}
\begin{equation*}
\inferrule*[right=$\Pi$-comp.]{
  y : \Pi(A, B) \t C(y) \ty \\
  f : (x : A)\, B(x) \t d(f) : C(\lambda(f)) \\
  x : A \t g(x) : B(x)}{
  \textsf{funsplit}(d\c \lambda(g)) = d(g) : C(\lambda(g))}
\end{equation*}

\pgph Observe that the assumption $f : (x : A)\,B(x) \t d(f) : C(\lambda(f))$ makes the
$\mathsf{funsplit}$ rules into higher-order inference rules, which as such are inexpressible
in the ``polymorphic'' formulation of type theory. Our task in the remainder of this section
will be to reformulate these rules in a first-order fashion. Our treatment is a
generalisation of that given by Martin-L\"of in his introduction to \cite{pml}, with the major
difference that we are working in the theory with \emph{intensional} identity types, as
opposed to the \emph{extensional} equality types of \cite{pml}.

\begin{Prop}(cf. \cite[p.\ 52]{Programming})\label{defapp}
In the presence of the rules $\Pi$-\textsc{form}, $\Pi$-\textsc{intro},
\mbox{$\Pi$-\textsc{elim}} and $\Pi$-\textsc{comp}, the rules $\Pi$-\textsc{app} and
$\Pi$-$\beta$ are definable.
\end{Prop}
\begin{proof}
Suppose that $m : \Pi(A, B)$ and $a : A$. We define a type $y : \Pi(A, B) \t
C(y) \ty$ by taking $C(y) \defeq B(a)$; and a term $f : (x : A)B(x) \t d(f) :
C(\lambda(f))$ by taking $d(f) \defeq f(a)$. Applying $\Pi$-elimination, we
define $\textsf{app}(m, a) \defeq \textsf{funsplit}(d, m) : B(a)$.
Moreover, when $m = \lambda(f)$ we have by $\Pi$-\textsc{comp} that
$\textsf{app}(\lambda(f), a) = d(f) = f(a)$, which gives us $\Pi$-$\beta$ as
required.
\end{proof}

\begin{Prop}(cf. \cite[p.\ 62]{Programming})\label{defeta}
In the presence of the  identity types and the  rules
$\Pi$-\textsc{form}, $\Pi$-\textsc{intro}, $\Pi$-\textsc{elim} and $\Pi$-\textsc{comp}, the
following rules are definable:
\begin{mathpar}
\inferrule*[right=$\Pi$-prop-$\eta$;]{
  m : \Pi(A, B)}{
  {\eta}(m) : \Id_{\Pi(A, B)}\big(m\c \lambda x.\, m \cdot x\big)}

\inferrule*[right=$\Pi$-prop-$\eta$-comp.]{
    x : A \t f(x) : B(x)}{
  {\eta}(\lambda(f)) = r(\lambda(f)) : \Id_{\Pi(A, B)}\big(\lambda(f)\c \lambda(f)\big)}
\end{mathpar}
\end{Prop}
\begin{proof}
  Given $y : \Pi(A, B)$, we define a type $C(y) \defeq \Id_{\Pi(A,
    B)}(y\c \lambda x.\, y \cdot x)$. In the case where $y =
  \lambda(f)$ for some $f : (x : A)B(x)$, we have $C(y) =
  {\Id_{\Pi(A, B)}(\lambda(f)\c \lambda x.\, \lambda(f) \cdot x)} =
  \Id_{\Pi(A, B)}(\lambda(f)\c \lambda x.\, f(x)) = \Id_{\Pi(A,
    B)}(\lambda(f)\c \lambda(f))$; consequently, we may define an element \mbox{$d(f) :
  C(\lambda(f))$} by $d(f) \defeq r(\lambda(f))$. Using
  $\Pi$-elimination, we define ${\eta}(m) \defeq \mathsf{funsplit}(d,
  m)$; and when $m = \lambda(f)$, we have by $\Pi$-\textsc{comp} that
  ${\eta}(\lambda(f)) = d(f) = r(\lambda(f))$ as required.
\end{proof}
\begin{Prop}\label{deffunsplit}
In the presence of identity types, the  rules $\Pi$-\textsc{form},
\mbox{$\Pi$-\textsc{intro}}, \mbox{$\Pi$-\textsc{app}} and $\Pi$-$\beta$, and the rules
$\Pi$-\textsc{prop}-$\eta$ and $\Pi$-\textsc{prop}-$\eta$-\textsc{comp} of Proposition
\ref{defeta}, the rules $\Pi$-\textsc{elim} and $\Pi$-\textsc{comp} are
definable.
\end{Prop}
\begin{proof}
We first recall that in the presence of identity types, $\Pi$-\textsc{form},
$\Pi$-\textsc{intro}, $\Pi$-\textsc{app} and $\Pi$-$\beta$, we may derive the following ``Leibnitz rules'', which assuming $A \ty$ and $x : A \t B(x) \ty$, say that
\begin{mathpar}
\inferrule*[right=$\Id$-subst;]{
  a_1, a_2 : A \\
  p : \Id_A(a_1, a_2) \\
  b_2 : B(a_2)
  }{
  \mathsf{subst}(p, b_2) : B(a_1)}

\inferrule*[right=$\Id$-subst-comp;]{
  a : A \\
  b : B(a)
  }{
  \mathsf{subst}(r(a), b) = b : B(a)}
\end{mathpar}
see \cite[p.\ 59]{Programming}, for example.

So, suppose given judgements $A \ty$, \mbox{$x : A \t B(x) \ty$} and $y :
\Pi(A, B) \t C(y) \ty$ and terms \mbox{$f : (x : A)B(x) \t d(f) :
C(\lambda(f))$} and $m : \Pi(A, B)$. We are required to define a term
$\mathsf{funsplit}(d, m) : C(m)$ satisfying $\mathsf{funsplit}(d,
\lambda(f)) = d(f)$. We begin by forming the term
\begin{equation*}
    T(d, m) \defeq d\big([x : A]\,m \cdot x\big) : C(\lambda x.\, m \cdot x)\text.
\end{equation*}
By $\Pi$-\textsc{prop}-$\eta$, we have a term $
  {\eta}(m) : \Id_{\Pi(A, B)}(m\c \lambda x.\, m \cdot x)$:
so by substituting $T(d, m)$ along ${\eta}(m)$ we obtain a term \mbox{$
    \mathsf{funsplit}(d, m) \defeq \mathsf{subst}\big({\eta}(m)\c T(d, m)\big) : C(m)$} as required.
Moreover, when $m = \lambda(f)$, we obtain from
    $\Pi$-$\beta$ that $T(d\c \lambda(f)) = d(f)$, and from $\Pi$-\textsc{prop}-$\eta$-\textsc{comp} that
    ${\eta}(\lambda(f)) = r(\lambda(f))$; and so from \mbox{$\Id$\textsc{-subst-comp}}, we deduce
    that
\begin{equation*}
    \mathsf{funsplit}(d\c \lambda(f)) = \mathsf{subst}\big(r(\lambda(f))\c d(f)\big) = d(f) : C(\lambda(f))
\end{equation*}
as required.
\end{proof}

\pgph \label{propequiv} Thus, in the presence of identity types, the
\textsf{funsplit}-formulation of dependent products is equivalent with the
\textsf{app}-formulation extended with the propositional $\eta$-rule. Note
carefully that this equivalence is a \emph{propositional}, rather than
\emph{definitional} one; which is to say that if we are given a
$\textsf{funsplit}$ term, to which we apply Propositions~\ref{defapp} and
\ref{defeta} to obtain terms $\textsf{app}$ and $\eta$, and then
Proposition~\ref{deffunsplit} to obtain a new term $\textsf{funsplit}'$, we
should not expect $\textsf{funsplit}(d, m) = \textsf{funsplit}'(d, m)$ to hold,
but rather only that
\begin{equation*}
\inferrule{
  y : \Pi(A, B) \t C(y) \ty \\
  f : (x : A)B(x) \t d(f) : C(\lambda(f)) \\
  m : \Pi(A, B)}{
  \psi(d, m) : \Id_{C(m)}\big(\textsf{funsplit}(d, m), \textsf{funsplit}'(d, m)\big)}
\end{equation*}
should be derivable. We may prove this by an application of
$\Pi$-\textsc{elim}.

\section{$\Pi$-application does not entail $\Pi$-elimination}
\pgph We saw in Proposition~\ref{defapp} that the \textsf{funsplit}-formulation
of dependent products subsumes the \textsf{app}-formulation; and the task of
this section is to show that the converse does not obtain. In the previous
section we were proving a positive derivability result, and so worked in a
minimal fragment of Martin-L\"of type theory in order to make our result as
strong as possible.  In this section, we are proving a negative derivability
result: and to make this as strong as possible, we work in full Martin-L\"of
type theory. So in addition to identity types and the
$\mathsf{app}$-formulation of dependent products we assume the presence of
dependent sums $\Sigma x :A.\, B(x)$, the unit type $\mathsf 1$, pairwise disjoint unions $A + B$, the empty type $\mathsf 0$, the $\mathsf W$-types, and the first universe $\mathsf U$. We refer to the
type theory with these constructors as $\cat{ML}_\textsf{app}$. Our main result
will be:

\begin{Thm}\label{main} Relative to the theory $\cat{ML}_\mathsf{app}$, the $\mathsf{funsplit}$ rules $\Pi$-\textsc{elim} and $\Pi$-\textsc{comp} are not
definable.
\end{Thm}
Now, if we could define $\Pi$-\textsc{elim} and $\Pi$-\textsc{comp} relative to
$\cat{ML}_\textsf{app}$, then by Proposition~\ref{defeta} we would also be able
to derive $\Pi$-\textsc{prop}-$\eta$ and
$\Pi$-\textsc{prop}-$\eta$-\textsc{comp}. Consequently, we may prove Theorem
\ref{main} by proving:

\paragraph*{{\mdseries\ref{main}'.}\ \ Theorem:}{\em
Relative to the theory $\cat{ML}_\mathsf{app}$, the rules $\Pi$-\textsc{prop}-$\eta$ and
$\Pi$-\textsc{prop}-$\eta$-\textsc{comp} of Proposition~\ref{defeta} are not definable.}
\vskip\baselineskip

\pgph Our method of proving Theorem \ref{main}' will be as follows. We first define the
following rules relative to the theory $\cat{ML}_\textsf{app}$:
\begin{mathpar}
\inferrule*[right=$\Pi'$-form;]{A \ty \\ x : A \t B(x) \ty}{\Pi'(A, B) \ty}

\inferrule*[right=$\Pi'$-abs;]{
x : A \t f(x) : B(x)}{\lambda'(f) : \Pi'(A, B)}

\inferrule*[right=$\Pi'$-app;]{
  m : \Pi(A, B) \\
  a : A}{
  \textsf{app}'(m, a) : B(a)}

\inferrule*[right=$\Pi'$-$\beta$.]{
  x : A \t f(x) : B(x) \\
  a : A}{
  {\textsf{app}'(\lambda'(f), a) = f(a) : B(a)}}
\end{mathpar}
We then show that the corresponding rule
$\Pi'$-\textsc{prop}-$\eta$ is not definable; and from this we deduce that the
rule $\Pi$-\textsc{prop}-$\eta$ cannot be definable either, since if it were
then by replacing each $\Pi$, $\lambda$ or $\textsf{app}$ in its derivation
with a $\Pi'$, $\lambda'$ or $\mathsf{app}'$, we would obtain a derivation of
$\Pi'$-\textsc{prop}-$\eta$, which would give a contradiction.

\newcommand{\inl}{{\text{\footnotesize{$\amalg$}}}_1}
\newcommand{\inr}{{\text{\footnotesize{$\amalg$}}}_2}

\pgph In order to define $\Pi'$, $\lambda'$ and $\mathsf{app}'$, we will need
to make use of \emph{disjoint union} types. Given types $A$ and $B$, their
disjoint union is the type $A + B$ with the following introduction and
elimination rules:
\begin{mathpar}
\inferrule*[right=$+$-intro$_1$;]{a : A}{\inl(a) : A + B}

\inferrule*[right=$+$-intro$_2$;]{b : B}{\inr(b) : A + B}

\inferrule*[right=$+$-elim{,}]{
  z : A + B \t C(z) \ty \\
  x : A \t f(x) : C(\inl(x))\\
  y : B \t g(y) : C(\inr(y))\\
  c : A + B
  }{\mathsf{case}(f, g, c) : C(c)}
\end{mathpar}
subject to the computation rules $\mathsf{case}\big(f, g, \inl(a)\big) = f(a)$ and
$\mathsf{case}\big(f, g, \inr(b)\big) = g(b)$. We use disjoint unions to define the $\Pi'$-rules as follows.
\begin{mathpar}
 \inferrule*[right=$\Pi'$-form;]{A \ty \\ x : A \t B(x) \ty}{\Pi'(A, B) \defeq \Pi(A, B) + \Pi(A, B) \ty}

\inferrule*[right=$\Pi'$-abs;]{
  x : A \t f(x) : B(x)}{\lambda'(f) \defeq \inl(\lambda(f)) : \Pi(A, B) + \Pi(A, B)}

\inferrule*[right=$\Pi'$-app,]{
  m : \Pi(A, B) + \Pi(A, B)\\
  a : A}{
  \textsf{app}'(m, a) \defeq \mathsf{case}(\textsf{app}(\thg, a)\c \textsf{app}(\thg, a)\c m) : B(a)}
\end{mathpar}
where we write $\textsf{app}(\thg, a)$ as an abbreviation for the term ${[x :
\Pi(A, B)]}\, \textsf{app}(x, a)$. To see that these definitions validate
$\Pi'$-$\beta$, we suppose that $f : (x : A)\, B(x)$ and $a : A$; then
by the first computation rule for disjoint unions and \textsc{$\Pi$-$\beta$} we have
that
\begin{equation*}
    \textsf{app}'(\lambda'(f), a) = \mathsf{case}\big(\textsf{app}(\thg, a)\c \textsf{app}(\thg, a)\c \inl(\lambda(f))\big) =
\textsf{app}(\lambda(f), a) = f(a)
\end{equation*}
as required.
%
%
%
%
%

\pgph It remains to show that with respect to the above definitions, the rule
\begin{equation*}
 \inferrule*[right=$\Pi'$-prop-$\eta$]{
  m : \Pi'(A, B)}{
  {\eta'}(m) : \Id_{\Pi'(A, B)}\big(m\c \lambda' x.\, \mathsf{app'}(m, x)\big)}
\end{equation*}
cannot be derived. So suppose that it could. Since for each judgement $x : A \t
f(x) : B(x)$ we have a term $\inr(\lambda(f)) : \Pi'(A, B)$, we would obtain
from this a derivation of
\begin{equation*}
 \inferrule{
  x : A \t f(x) : B(x)}{
  {\eta'}(\inr(\lambda(f))) : \Id_{\Pi'(A, B)}\big(\inr(\lambda(f))\c \lambda' x.\, \mathsf{app'}(\inr(\lambda(f)), x)\big)\text.}
\end{equation*}
But now by the definition of $\textsf{app}'$ we have that
\begin{equation*}
    \textsf{app}'(\inr(\lambda(f)),x) = \mathsf{case}(\textsf{app}(\thg, x)\c
\textsf{app}(\thg, x)\c \inr(\lambda(f))) = \textsf{app}(\lambda(f), x) =
f(x)\text;
\end{equation*}
and hence $\lambda' x.\, \mathsf{app'}(\inr(\lambda(f)), x) = \lambda' x.\, f(x) = \inl(\lambda(f))$, so that we may view the above derivation as a
derivation of
\begin{equation}\label{specialtype}\tag{$\star$}
 \inferrule{
   A \ty \\
  x : A \t B(x) \ty \\
  x : A \t f(x) : B(x)}
  {\eta'(\inr(\lambda(f))) : \Id_{\Pi(A, B) + \Pi(A, B)}\big(\inr(\lambda(f))\c \inl(\lambda(f))\big)\text.}
\end{equation}
To complete the proof, it suffices to show that no such derivation can exist. The key to doing
so will be the following \emph{disjointness rule}:
\begin{equation}\label{contra}\tag{$\dagger$}
 \inferrule{
   C \ty \\
  c : C\\ p : \Id_{C + C}\big(\inr(c)\c \inl(c)\big)}
  {\theta(c, p) : \mathsf 0}\text,
\end{equation}
where $\mathsf 0$ is the empty type. If we can prove that this holds relative to
$\cat{ML}_\textsf{app}$, then we will be able to deduce the underivability of
\eqref{specialtype}. Indeed, suppose that \eqref{specialtype} holds. Then from this and
\eqref{contra} we can derive the following rule:
\begin{equation*}
 \inferrule{
  x : A \t f(x) : B(x)}
  {\theta_{\Pi(A, B)}\big(\lambda(f)\c \eta'(\inr(\lambda(f)))\big) : \mathsf 0}\text;
\end{equation*}
and by instantiating this derivation at some particular $A$, $B$ and $f$---a suitable choice
being $A
\defeq \mathsf 1$, $B
\defeq \mathsf 1$ and $f
\defeq [x : \mathsf 1]\, x$---we obtain a global
element of $\mathsf 0$. But this is impossible, because $\cat{ML}_\textsf{app}$
is known to be \emph{consistent}, in the sense that $\mathsf 0$ has \emph{no}
global elements. An easy way of seeing that this is the case is by exhibiting a consistent model
for $\cat{ML}_\textsf{app}$ using the sets in our meta-theory. We interpret
types as sets; dependent sums and products as indexed sums and products;
identity types as meta-theoretic equality; the terminal type as a one-element
set and the empty type as the empty set. The interpretation of $\mathsf
W$-types and the first universe is a little more complex, and depends upon the
existence of inductive datatypes in our meta-theory, but is essentially
unproblematic.

\pgph \label{deferlemma} All that remains to complete the proof of Theorem \ref{main}' is to
show that the disjointness rule \eqref{contra} is derivable in $\cat{ML}_\textsf{app}$. This
follows by a standard argument (cf.~\cite[p.~86]{Programming}). Recall that one of the type
constructors in $\cat{ML}_\textsf{app}$ was that for the \emph{universe type}~\cite[Chapter
14]{Programming}. This is a type $\mathsf U$ containing ``codes'' for each of the other type
formers of $\cat{ML}_\textsf{app}$. In particular, we have rules
\begin{equation*}
\vcenter{\hbox{$\inferrule*[right=$\mathsf U$-intro$_1$]{\quad\qquad}{\hat {\mathsf 0} :
\mathsf U}$}}
\qquad \text{and} \qquad
 \vcenter{\hbox{$\inferrule*[right=$\mathsf U$-intro$_2$]{\quad\qquad}{\hat {\mathsf
1} : \mathsf U}$}}
\end{equation*}
introducing codes for the empty type and the terminal type. Recall also that
$\mathsf U$ comes equipped with a \emph{decoding function} $\mathsf D$ which is
given by an indexed family of types
\begin{equation*}
    x : \mathsf U \t \mathsf D(x) \ty
\end{equation*}
together with computation rules which determine the value of $\mathsf D$ on
canonical elements of $\mathsf U$. In particular, we have rules
\begin{equation*}
\vcenter{\hbox{$\inferrule*[right=$\mathsf U$-comp$_1$]{\quad\qquad}{\mathsf D(\hat {\mathsf 0}) = \mathsf 0 \ty}$}}
\qquad \text{and} \qquad
\vcenter{\hbox{$\inferrule*[right=$\mathsf U$-comp$_2$\text.]{\quad\qquad}{\mathsf D(\hat {\mathsf 1}) = \mathsf 1 \ty}$}}
\end{equation*}
So suppose now that $C \ty$, $c : C$ and $p : \Id_{C + C}\big(\inr(c)\c
\inl(c)\big)$ as in the premisses of~\eqref{contra}. We are required to derive
an element of $\mathsf 0$. We begin by defining functions
\begin{equation*}
    x : C \t f(x) \defeq \hat{\mathsf 0} : \mathsf U \qquad \text{and} \qquad x : C \t g(x) \defeq \hat{\mathsf 1} : \mathsf U\text.
\end{equation*}
Applying $+$-elimination to these we obtain a function $\textsf{case}(f, g, \thg) \colon C + C
\to \mathsf U$; and using the decoding function $\mathsf D$ on this we obtain a family
\begin{equation*}
    z : C + C \t T(z) \defeq \mathsf D\big(\textsf{case}(f, g, z)\big) \ty\text.
\end{equation*}
Now from the rule \textsc{$\Id$-subst} defined in
Proposition~\ref{deffunsplit}, together with the given proof $p : \Id_{C +
C}\big(\inr(c)\c \inl(c)\big)$ we obtain the term
\begin{equation*}
    x : T(\inl(c)) \t \mathsf{subst}(p, x) : T(\inr(c))\text.
\end{equation*}
But we have that $T(\inl(c)) = \mathsf D (f(c))  = \mathsf D(\hat{\mathsf 0}) = \mathsf 0$ and
that $T(\inr(c)) = \mathsf D(g(c)) = \mathsf D(\hat{\mathsf 1}) = \mathsf 1$, so that we may
view this as a function $ x : \mathsf 1  \t \mathsf{subst}(p, x) : \mathsf 0$. In particular,
by evaluating this function at the canonical element $\mathord \star : \mathsf 1$ we obtain an
element  $\mathsf{subst}(p, \mathord \star) : \mathsf 0$ as required. This completes the proof
of Theorem \ref{main}'.

\section{Function extensionality}
\pgph \label{funextdef} In this final section, we investigate the principle of \emph{function
extensionality} in Martin-L\"of type theory, which asserts that two elements of a dependent
product type which are pointwise propositionally equal, are themselves propositionally equal.
Explicitly, it is given by the following two inference rules:
\begin{mathpar}
\inferrule*[right=$\Pi$-ext{,}]{
  m\c n : \Pi(A, B)\\
  k : \Pi x: A.\, \Id_{B(x)}(m \cdot x,\, n \cdot x)}
  {\textsf{ext}(m, n, k) : \Id_{\Pi(A, B)}(m, n)}

\inferrule*[right=$\Pi$-ext-comp{,}]{
  f : (x : A)\, B(x)}
  {\textsf{ext}\big(\lambda(f), \lambda(f), \lambda(rf)\big) = r(\lambda(f)) : \Id_{\Pi(A, B)}(\lambda(f), \lambda(f))}
\end{mathpar}
where we write $rf$ as an abbreviation for the term $[x:A]\, rf(x)$. These
rules were considered first by Turner in~\cite{DT} and then more extensively by
Hofmann~\cite{Hof}. If one views Martin-L\"of type theory as a
\emph{computational} system, in which terms are thought of as algorithms---an
idea made precise in \cite[Appendix B]{Programming}, for example---then these
rules are hard to justify, since two extensionally equal functions can have
quite different algorithmic content. From a proof-theoretic perspective they
are also problematic, since they destroy one of the more pleasant properties of
Martin-L\"of type theory, namely that in the syntactic model, every global
element of a closed type is definitionally equal to a canonical
element.\footnote{Though a construction of Altenkirch~\cite{alt} shows that
there are models validating both function extensionality and the canonicity
property.} On the other hand, they do not break strong normalisation, so that
if we view type theory merely as a \emph{computable} system---one for which the
correctness of derivations is decidable---then their addition is unproblematic,
and in fact produces a system which is closer to the ``everyday'' mathematics
described by extensional type theory. Indeed, Hofmann~\cite{Hof2} shows that
augmenting intensional type theory with function extensionality and the
principle of \emph{uniqueness of identity proofs}---which asserts that any two
proof-terms $p, q : \Id_A(a, b)$ are themselves propositionally equal---yields
a system which, whilst still decidable, is propositionally equivalent (in the
sense of~\S\ref{propequiv}) to extensional type theory.

The author's motivations for studying the principle of function extensionality
are somewhat different from those of~\cite{Hof2}; they are informed by his
investigations~\cite{g1} into two-dimensional semantics for type theory. In
this semantics, dependent product formation is required to be a (suitably weak)
two-dimensional right adjoint to reindexing; and in order for the semantics to
be complete, we must verify that the syntactic model has this property---which
requires the imposition of some form of function extensionality. However,
whilst preparing~\cite{g1}, it became apparent to the author that the usual
formulation of function extensionality is insufficient for this purpose because
it fails to verify a number of necessary propositional equalities between
identity proofs. In the setting of~\cite{Hof2}, the existence of these
propositional equalities is assured by the principle of uniqueness of identity
proofs; yet for a higher-dimensional semantics it is crucial that we do
\emph{not} have uniqueness of identity proofs, whose imposition would allow
only trivial, posetal, models to be formed. Thus it is of interest to determine
how function extensionality should correctly be formulated when we do not have
uniqueness of identity proofs; and it is this that we shall now do. We work in
the fragment of type theory given by the identity types and the
\textsf{app}-formulation of dependent products. In order to minimize clutter,
we also allow ourselves the notational convenience of writing function
application $f(x)$ simply as $fx$, and $\lambda$-abstraction $\lambda(f)$
simply as $\lambda f$. We begin by recording some useful consequences of
function extensionality.

\begin{Prop}
In the presence of $\Pi$-\textsc{ext} and $\Pi$-\textsc{ext}-\textsc{comp}, the rules
$\Pi$-\textsc{prop}-$\eta$ and $\Pi$-\textsc{prop}-$\eta$-\textsc{comp} of
Proposition~\ref{defeta} are definable.
\end{Prop}
\begin{proof}
Given $m : \Pi(A, B)$, we must exhibit a term
    ${\eta}(m) :
\Id_{\Pi(A, B)}(m\c \lambda x.\, m \cdot x)$. So we define $n = \lambda x.\, m \cdot x : \Pi(A, B)$; and by $\Pi$-$\beta$ have that $n \cdot x =
(\lambda x.\, m \cdot x) \cdot x = m \cdot x$ whenever $x : A$. We may now
define $ \eta(m)
\defeq \textsf{ext}\big(m\c n\c \lambda x.\, r(m \cdot x)\big) $; and moreover,
when $m = \lambda f$ for some $f : (x : A)\, B(x)$, the $\beta$-rule implies
that $m = n$, so that $\eta(\lambda f) = \textsf{ext}(\lambda f, \lambda f,
\lambda(rf)) = r(\lambda f)$ as required.
\end{proof}
\begin{Prop}\label{defxi}
In the presence of $\Pi$-\textsc{ext} and $\Pi$-\textsc{ext}-\textsc{comp}, the following
\emph{propositional $\xi$-rules} are definable:
\begin{mathpar}
\inferrule*[right=$\Pi$-prop-$\xi${,}]{
  f\c g : (x : A)\, B(x)\\
  p : (x : A)\, \Id_{B(x)}(fx, gx)}
  {\xi(f, g, p) : \Id_{\Pi(A, B)}(\lambda f, \lambda g)}

\inferrule*[right=$\Pi$-prop-$\xi$-comp{.}]{
  f : (x : A)\, B(x)}
  {\xi(f, f, rf) = r(\lambda f) : \Id_{\Pi(A, B)}(\lambda f, \lambda f)}
  \end{mathpar}
\end{Prop}
\begin{proof}
Given $f$, $g$ and $p$ as in the hypotheses of $\Pi$-\textsc{prop}-$\xi$, we
consider $m = \lambda f$ and $n = \lambda g$ in $\Pi(A, B)$. By the
$\beta$-rule, we may view $p$ as a term
\begin{equation*}
    x : A \t p(x) :
\Id_{B(x)}(m \cdot x\c n \cdot x)\text;
\end{equation*}
and hence may define $\xi(f, g, p) = \mathsf{ext}(\lambda f, \lambda g,
\lambda p)$. Moreover, we have that $\xi(f, f, rf) = \textsf{ext}(\lambda f,
\lambda f, \lambda(rf)) = r(\lambda f)$ as required.
\end{proof}

In fact, we have a converse to the previous two propositions:

\begin{Prop}\label{referens}
In the presence of the rules $\Pi$-\textsc{prop}-$\eta$ and
$\Pi$-\textsc{prop}-$\eta$-\textsc{comp} of Proposition~\ref{defeta} and the rules
$\Pi$-\textsc{prop}-$\xi$ and \mbox{$\Pi$-\textsc{prop}-$\xi$-\textsc{comp}} of
Proposition~\ref{defxi}, the function extensionality rules $\Pi$-\textsc{ext} and
$\Pi$-\textsc{ext}-\textsc{comp} are definable.
\end{Prop}
\begin{proof}
Recall from \cite{HS} that, in the presence of dependent products, the identity
types admit an operation which one may think of as either \emph{transitivity}
or \emph{composition}:
\begin{mathpar}
\inferrule*[right=$\Id$-trans{,}]{
  p : \Id_A(a_1, a_2) \\
  q : \Id_A(a_2, a_3)
  }{
  q \circ p : \Id_A(a_1, a_3)}

\inferrule*[right=$\Id$-trans-comp;]{
  p : \Id_A(a_1, a_2)
  }{
  p \circ r(a_1) = p : \Id_A(a_1, a_2)}
\end{mathpar}
and also an operation which one may think of as either \emph{symmetry} or
\emph{inverse}:
\begin{mathpar}
\inferrule*[right=$\Id$-symm{,}]{
  p : \Id_A(a_1, a_2)
  }{
  p^{-1} : \Id_A(a_2, a_1)}

\inferrule*[right=$\Id$-symm-comp.]{
  a : A
  }{
  r(a)^{-1} = r(a) : \Id_A(a, a)}
\end{mathpar}
Now suppose that we are given terms $m$, $n$ and $k$ as in the hypotheses of
$\Pi$\textsc{-ext}. We begin by defining terms
\begin{align*}
    f &\defeq [x : A]\, m \cdot x : (x : A)\, B(x)\\
    g &\defeq [x : A]\, n \cdot x : (x : A)\, B(x)\\
    \text{and} \qquad p & \defeq [x : A]\, k \cdot x : (x : A)\, \Id_{B(x)}(fx, gx)\text.
\end{align*}
Observe that the third of these is well-typed by virtue of the first two.
Applying the propositional $\xi$-rule, we obtain a term
 \begin{equation*}
    \xi(f, g, p) : \Id_{\Pi(A, B)}(\lambda f, \lambda g) = \Id_{\Pi(A, B)}\big(\lambda x.\, m \cdot x\c \lambda x.\, n \cdot x\big)\text.
 \end{equation*}
But from the propositional $\eta$-rule and $\Id$-symmetry rule, we have terms
 \begin{gather*}
    {\eta}(m)  : \Id_{\Pi(A, B)}(m\c \lambda x.\, m \cdot x) \quad \text{and} \quad
    {\eta}(n)^{-1} : \Id_{\Pi(A, B)}(\lambda x.\, n \cdot x\c n)
 \end{gather*}
and now can define
    $\textsf{ext}(m,n,p) \defeq \eta(n)^{-1} \circ \big(\xi(f, g, p) \circ \eta(m)\big) : \Id_{\Pi(A,
    B)}(m,n)$.
In the case where $m = n = \lambda h$ and $p = \lambda(rh)$ we have by the
$\beta$-rule that $f = g = h$, and so we may calculate that
\begin{align*}
    \textsf{ext}(\lambda h, \lambda h, \lambda(rh)) &= \eta(\lambda h)^{-1} \circ \big(\xi(h, h, rh) \circ \eta(\lambda h)\big)\\
    & = r(\lambda h)^{-1} \circ \big(r(\lambda h) \circ r(\lambda h)\big) = r(\lambda h)
\end{align*}
as required.
\end{proof}
Thus relative to the theory with identity types plus the
$\textsf{app}$-formulation of dependent products, the function extensionality
principle is equivalent\footnote{Again, in a propositional, rather than
definitional, sense.} with the conjunction of the propositional $\eta$-\ and
propositional $\xi$-rules; and relative to the theory with identity types plus
the $\textsf{funsplit}$ formulation of dependent products, function
extensionality is equivalent with the propositional $\xi$-rule.

%
%
%
%

\pgph We now wish to describe the inadequacies of function extensionality in
the absence of uniqueness of identity proofs. These arise from its failure to
continue a characteristic trend in intensional type theory, namely that nearly
every statement that one may think should hold about the identity types, does
hold. For instance, in the proof of Proposition~\ref{referens}, we saw that the
identity types $\Id_A$ come equipped with operations which we called
\emph{composition} and \emph{inverse}. We would hope for this composition to be
associative and unital, and for the inverse operation to really provide
compositional inverses; and a straightforward application of $\Id$-elimination
shows this to be the case, at least when we interpret associativity, unitality
and invertibility in an ``up-to-propositional-equality'' sense. Similarly, each
judgement $x : A \t f(x) : B(x)$ induces a judgement $x, y : A\c p : \Id_A(x,
y) \t \tilde f(p) : \Id_{B(x)}(fx, fy)$ which we would expect to be suitably
``functorial'' in $p$: and again, an application of $\Id$-elimination confirms
this, providing us with canonical propositional equalities between $\tilde f(q
\circ p)$ and $\tilde f(q) \circ \tilde f(p)$. However, when it comes to
function extensionality, there are a number of statements which intuitively
should be true but which seem to be impossible to prove. Here are two typical
examples.

\begin{Exs}\label{section56}
\begin{enumerate}
\item Using $\Id$-elimination we can derive a rule
\begin{equation*}
\inferrule{
  m\c n : \Pi(A, B)\\ p : \Id_{\Pi(A, B)}(m, n) \\ a : A}
  {p \ast a  : \Id_{B(a)}(m \cdot a\c n \cdot a)}
\end{equation*}
satisfying $r(m) \ast a = r(m \cdot a)$, which expresses that any two
propositionally equal elements of a $\Pi$-type are pointwise
propositionally equal. We would expect that for $k : \Pi x : A.\,
\Id_{B(x)}(m \cdot x, n \cdot x)$ and $a : A$, we should have $k \cdot a$
propositionally equal to $\textsf{ext}(m, n, k) \ast a$; yet this seems
impossible to prove. \vskip\baselineskip


\item Suppose given terms $\ell, m, n : \Pi(A, B)$ and proofs \mbox{$f : (x
    : A)\, \Id_{B(x)}(\ell \cdot x,\, m \cdot x)$} and $g :  (x : A)\,
    \Id_{B(x)}(m \cdot x,\, n \cdot x)$. Let us write $g \circ f$ for the
    term $[x : A]\, gx \circ fx$. It now seems to be impossible to
    verify a propositional equality between the elements
\begin{equation*}
    \textsf{ext}(m, n, \lambda g) \circ \textsf{ext}(\ell, m, \lambda f)
\quad \text{and} \quad
    \textsf{ext}(\ell, n, \lambda(g \circ f))
\end{equation*}
of $\Id_{\Pi(A, B)}(\ell, n)$.
\end{enumerate}
\end{Exs}

\pgph \label{newrule} The reason that we encounter these problems is
essentially the following. We would like to construct the desired propositional
equalities by eliminating over the type \mbox{$u, v : \Pi(A, B) \t \Pi x: A.\,
\Id_{B(x)}(u \cdot x,\, v \cdot x)$}. To do this we need an elimination rule
that we do not have, one which says that this type is generated by elements of
the form $(\lambda f, \lambda f, \lambda(rf))$. In light of this, we propose
that function extensionality should be replaced with just such an elimination
rule. We consider the following two rules:
\begin{mathpar}
\inferrule*[right=$\Pi$-$\Id$-elim{,}]{
  u, v : \Pi(A, B)\c
  w : \Pi x: A.\, \Id_{B(x)}(u \cdot x,\, v \cdot x) \t C(u, v, w) \ty \\ f : (x : A)\,B(x) \t d(f) : C(\lambda f, \lambda f, \lambda(rf))\\
  m, n : \Pi(A, B) \\ k : \Pi x: A.\, \Id_{B(x)}(m \cdot x,\, n \cdot x)
  }
  {L(d, m, n, k) : C(m, n, k)}

\inferrule*[right=$\Pi$-$\Id$-comp{.}]{
  u, v : \Pi(A, B)\c
  w : \Pi x: A.\, \Id_{B(x)}(u \cdot x,\, v \cdot x) \t C(u, v, w) \ty \\ f : (x : A)\,B(x) \t d(f) : C(\lambda f, \lambda f, \lambda(rf))\\
  h : (x : A)B(x)
  }
  {L(d, \lambda h, \lambda h, \lambda(rh)) = d(h) : C(\lambda h, \lambda h, \lambda(rh))}
\end{mathpar}
Observe that these two rules are once again higher-order inference rules. We
will return to this point in \S \ref{discretesection} below. Let us first show
that these rules entail function extensionality.
\begin{Prop}
In the presence of identity types, the $\mathsf{app}$-formulation of $\Pi$-types and the rules
\textsc{$\Pi$-$\Id$-elim} and \textsc{$\Pi$-$\Id$-comp} of \S \ref{newrule}, it is possible to
define the function extensionality rules \textsc{$\Pi$-ext} and \textsc{$\Pi$-ext-comp}.
\end{Prop}
\begin{proof}
For each $u, v : \Pi(A, B)$ and $w : \Pi x: A.\, \Id_{B(x)}(u \cdot x,\, v
\cdot x)$ we define a type $C(u, v, w) \defeq \Id_{\Pi(A, B)}(u, v)$; and now for
each $f : (x : A)\, B(x)$, we define an element \mbox{$d(f) \defeq
r(\lambda f) : C(\lambda f, \lambda f, \lambda(rf))$}. Applying
$\Pi$-$\Id$-elimination, we obtain the judgement
\begin{equation*}
  \inferrule{
  m\c n : \Pi(A, B)\\
  k : \Pi x: A.\, \Id_{B(x)}(m \cdot x,\, n \cdot x)}
  {\textsf{ext}(m, n, k) := L(d, m, n, k) : \Id_{\Pi(A, B)}(m, n)}\text;
\end{equation*}
and calculate that $\mathsf{ext}(\lambda f, \lambda f, \lambda(rf)) = L(d,
\lambda f, \lambda f, \lambda(rf)) = d(f) = r(\lambda f)$ as required.
\end{proof}
Thus the $\Pi$-$\Id$-elimination rules are stronger than the function
extensionality rules, and we conjecture that they are \emph{strictly} stronger.
To prove this would require either finding a new model of type theory that
supports function extensionality but not $\Pi$-$\Id$-elimination---new, because
every semantic model of which the author is aware that supports the former,
also supports the latter---or giving a syntactic proof along the lines of that
given for Theorem~\ref{main}. In both cases, the author's efforts have proven
fruitless. Setting aside this issue, let us now show that the
$\Pi$-$\Id$-elimination rules allow us to give positive answers to the question
posed in Examples \ref{section56}.
\begin{Prop}\label{specialrules}
In the presence of identity types, the $\mathsf{app}$-formulation of $\Pi$-types and the rules
\textsc{$\Pi$-$\Id$-elim} and \textsc{$\Pi$-$\Id$-comp} of \S \ref{newrule}, the following
rules are definable:
\begin{mathpar}
\inferrule*[right=$\Pi$-ext-app{,}]{
  m\c n : \Pi(A, B)\\
  k : \Pi x: A.\, \Id_{B(x)}(m \cdot x,\, n \cdot x) \\ a : A}
  {\mu(m, n, k, a) : \Id_{\Id_{B(a)}(m \cdot a,\, n \cdot a)}\big(\mathsf{ext}(m, n, k) \ast a\c k \cdot a\big)}

\inferrule*[right=$\Pi$-ext-app-comp{\text,}]{
  f : (x : A)\, B(x)}
  {\mu(\lambda f, \lambda f, \lambda(rf), a) = r(rfa) : \Id_{\Id_{B(a)}(fa, fa)}\big(rfa\c rfa\big)}
\end{mathpar}
where $\ast$ is the operation defined in Examples \ref{section56}(1).
\end{Prop}
\begin{proof}
For each $u, v : \Pi(A, B)$ and $w : \Pi x: A.\, \Id_{B(x)}(u \cdot x,\, v \cdot x)$ we define
a type $C(u, v, w) \defeq \Pi x : A.\, \Id_{\Id_{B(x)}(u \cdot x, v \cdot
x)}\big(\mathsf{ext}(u, v, w) \ast x\c w \cdot x\big)$. Now for $f : (x : A)\, B(x)$, we
calculate that
\begin{align*}
    C(\lambda f, \lambda f, \lambda(rf)) &= \Pi x : A.\, \Id_{\Id_{B(x)}(fx,
fx)}\big(\mathsf{ext}(\lambda f, \lambda f, \lambda(rf)) \ast
x\c rfx\big)\\
&= \Pi x : A.\, \Id_{\Id_{B(x)}(fx, fx)}\big(r(\lambda f) \ast
x\c rfx\big)\\
&= \Pi x : A.\, \Id_{\Id_{B(x)}(fx, fx)}\big(rfx\c rfx\big)
\end{align*}
so that we may define $d(f) := \lambda x.\, r(rfx) : C(\lambda f,
\lambda f, \lambda(rf))$. An application of $\Pi$-$\Id$-elimination now yields
the judgement $\Pi$-\textsc{ext-app} by taking
\begin{equation*}
\mu(m, n, k, a)
\defeq L(d, m, n, k) \cdot a : \Id_{\Id_{B(a)}(m \cdot a, n \cdot a)}\big(\mathsf{ext}(m, n, k) \ast a\c k
\cdot a\big)\text.
\end{equation*} Finally, we compute that $\mu(\lambda f, \lambda f, \lambda(rf), a) = \lambda
x.\, r(rfx) \cdot a = r(rfa)$ as required.
\end{proof}

\begin{Prop}
In the presence of identity types, the $\mathsf{app}$-formulation of $\Pi$-types and the rules
\textsc{$\Pi$-$\Id$-elim} and \textsc{$\Pi$-$\Id$-comp} of \S \ref{newrule}, the following
rule is definable:
\begin{equation*}
\inferrule{
    \ell\c m\c n : \Pi(A, B)\\
    f : (x : A)\, \Id_{B(x)}(\ell \cdot x, m \cdot x)\\ g : (x : A)\, \Id_{B(x)}(m \cdot x, n \cdot x)}
  {\nu(f, g) : \Id_{\Id_{\Pi(A, B)}(\ell, n)}\big(\mathsf{ext}(m, n, \lambda g) \circ \mathsf{ext}(\ell, m, \lambda f)\c \mathsf{ext}(\ell, n, \lambda(g \circ f))\big)}
\end{equation*}
\end{Prop}
\begin{proof}
It suffices to derive the rule:
\begin{equation*}
\inferrule{
    \ell\c m\c n : \Pi(A, B)\\
    j : \Pi x : A.\, \Id_{B(x)}(\ell \cdot x, m \cdot x)\\ k : \Pi x : A.\, \Id_{B(x)}(m \cdot x, n \cdot x)}
  {\nu'(j, k) : \Id_{\Id_{\Pi(A, B)}(\ell, n)}\big(\mathsf{ext}(m, n, k) \circ \mathsf{ext}(\ell, m, j)\c \mathsf{ext}(\ell, n, \lambda x.\, k \cdot x \circ j \cdot x)\big)}
\end{equation*}
since the required result then follows by taking $j := \lambda f$ and $k := \lambda g$. But
applying the rule of \mbox{$\Pi$-$\Id$-elimination} on $k$, it suffices to derive this rule in the case where $m
= n = \lambda h$ and $k = \lambda(rh)$; which is to show that
\begin{equation*}
\inferrule{
    \ell : \Pi(A, B)\\ h : (x : A)\, B(x) \\ j : \Pi x : A.\, \Id_{B(x)}(\ell \cdot x, hx)}
  {\nu'(j, \lambda(rh)) : \Id_{\Id_{\Pi(A, B)}(\lambda h, n)}\big(
  \mathsf{ext}(\lambda h\c \lambda h\c \lambda(rh)) \circ \mathsf{ext}(\ell\c \lambda h\c j)
  \c\\
  \quad \qquad \qquad \qquad  \quad \qquad \qquad \qquad \qquad
  \mathsf{ext}(\ell\c \lambda h\c \lambda x.\, r(hx) \circ j \cdot x)
  \big)}
\end{equation*}
is derivable. But we have that $r(hx) \circ j \cdot x = j \cdot x$ and that
$\mathsf{ext}(\lambda h, \lambda h, \lambda(rh)) = r(\lambda h)$ so that
$\mathsf{ext}(\lambda h, \lambda h, \lambda(rh)) \circ \mathsf{ext}(\ell,
\lambda h, j) = \mathsf{ext}(\ell, \lambda h, j)$: so that it suffices to
show that
\begin{equation*}
\inferrule{
    \ell : \Pi(A, B)\\ h : (x : A)\, B(x) \\ j : \Pi x : A.\, \Id_{B(x)}(\ell \cdot x, hx)}
  {\nu'(j, \lambda(rh)) : \Id_{\Id_{\Pi(A, B)}(\lambda h, n)}\big(\mathsf{ext}(\ell\c \lambda h\c j)\c
  \mathsf{ext}(\ell\c \lambda h\c \lambda x.\, j \cdot x)
  \big)}
\end{equation*}
is derivable. Now, using the propositional $\eta$-rule, we can derive a term $\eta(j)$
witnessing the propositional equality of $j$ and $\lambda x.\, j \cdot x$; and we will be
done if we can lift this to a propositional equality between $\mathsf{ext}(\ell, \lambda h,
j)$ and $\mathsf{ext}(\ell, \lambda h, \lambda x.\, j \cdot x)$. But we may do this using
the following rule:
\begin{equation*}
\inferrule{
    a, b : \Pi(A, B) \\ c, d : \Pi x : A.\, \Id_{B(x)}(a \cdot x, b \cdot x) \\ p : \Id_{\Pi x : A.\, \Id_{B(x)}(a \cdot x, b \cdot x)}(c, d)}
  {\widetilde{\mathsf{ext}}(p) : \Id_{\Id_{\Pi(A, B)}}\big(\mathsf{ext}(a, b, c)\c \mathsf{ext}(a, b, d)\big)\text,}
\end{equation*}
which is derivable by $\Id$-elimination on $p$.
\end{proof}

In Section 4, we saw that the higher-order formulation of $\Pi$-types can be
restated in a first-order manner; and the final result of this paper will do
something similar for the $\Pi$-$\Id$-elimination rule.

\begin{Prop}\label{discretesection}
In the presence of the identity types; the $\mathsf{app}$-formulation of $\Pi$-types; the
function extensionality rules \textsc{$\Pi$-ext} and \textsc{$\Pi$-ext-comp}; and the rules
\textsc{$\Pi$-ext-app} and \textsc{$\Pi$-ext-app-comp} of Proposition \ref{specialrules}, we
can define the rules \textsc{$\Pi$-$\Id$-elim} and \textsc{$\Pi$-$\Id$-comp} of \S
\ref{newrule}.
\end{Prop}
\begin{proof}
Suppose that we are given terms
\begin{equation*}
\inferrule{
  u, v : \Pi(A, B)\c
  w : \Pi x: A.\, \Id_{B(x)}(u \cdot x,\, v \cdot x) \t C(u, v, w) \ty \\ f : (x : A)B(x) \t d(f) : C(\lambda f, \lambda f, \lambda(rf))\\
  m, n : \Pi(A, B) \\ k : \Pi x: A.\, \Id_{B(x)}(m \cdot x,\, n \cdot x)
  }
  {}
\end{equation*}
as in the premisses of \textsc{$\Pi$-$\Id$-elim}. We must find an element $L(d,
m, n, k) : C(m, n, k)$. We will employ much the same method as we did in the
proof of Proposition \ref{deffunsplit}, though the details will be a little
more complicated. We begin by constructing the element \mbox{$d([x : A]\, m
\cdot x) : C(\lambda x.\, m \cdot x\c \lambda x.\, m \cdot x\c \lambda x.\, r(m
\cdot x))$}; and the remainder of the proof will involve applying various
substitutions to this element until we obtain the required element of $C(m, n,
k)$. The key result we need is the following lemma.

\begin{Lems}
We may define a rule:
\begin{equation}\tag{$\star$}\label{tagged}
\inferrule{
  u, v : \Pi(A, B)\\
  p : \Id_{\Pi(A, B)}(u, v)\\
  c : C(\lambda x.\, u \cdot x\c \lambda x.\, u \cdot x\c \lambda x.\, r(u \cdot x))
  }
  {\phi(p, c) : C(u, v, \lambda x.\, p \ast x)}
\end{equation}
satisfying $\phi\big(r(\lambda f)\c c\big) = c$.
\end{Lems}
Before proving this, let us see how it allows us to derive the required element of $C(m, n, k)$. Using
function extensionality we can form $\mathsf{ext}(m, n, k) : \Id_{\Pi(A, B)}(m, n)$; and so by
applying $\phi$ to this and $d([x : A]\, m \cdot x)$ can obtain an element
\begin{equation*}
    b(m, n, k) \defeq \phi\big(\textsf{ext}(m, n, k)\c d([x : A]\, m \cdot x)\big) : C(m, n, \lambda
    x.\, \textsf{ext}(m, n, k) \ast x)\text.
\end{equation*}
We now make use of the rule $\Pi$-\textsc{ext-app} of Proposition \ref{specialrules}, which
provides us with a term
\begin{equation*}
    x : A \t \mu(m, n, k, x) : \Id_{\Id_{B(x)}(m \cdot x,\, n \cdot x)}\big(\mathsf{ext}(m, n, k) \ast x\c k \cdot
    x\big)\text;
\end{equation*}
applying function extensionality to which yields a term
\begin{multline*}
    p(m, n, k) \defeq \mathsf{ext}(\lambda x.\, \mathsf{ext}(m, n, k) \ast x\c k\c \lambda x.\,\mu(m, n, k, x)\big) \\: \Id_{\Pi x: A.\, \Id_{B(x)}(m \cdot x,\, n \cdot
    x)}(\lambda x.\, \mathsf{ext}(m, n, k) \ast x\c k)\text.
\end{multline*}
The final step is to use the Leibnitz rule defined in the proof of Proposition
\ref{deffunsplit} to form the required term $L(d, m, n, k) \defeq \mathsf{subst}(p(m, n, k),
b(m, n, k)) : C(m, n, k)$.
 We are also required to show that $L(d, \lambda f, \lambda f,
\lambda(rf)) = d(f)$. For this, we first note that $b(\lambda f, \lambda f, \lambda(rf)) =
\phi(r(\lambda f), d(f)) = d(f) : C(\lambda f, \lambda f, \lambda(rf))$. Next we observe
that $\mu(\lambda f, \lambda f, \lambda(rf), x) = r(rfx)$ so that we have
\begin{align*}
p(\lambda f, \lambda f, \lambda(rf)) &= \mathsf{ext}\big(\lambda x.\, r(\lambda f) \ast x\c
\lambda(rf)\c \lambda x.\, r(rfx)\big)\\
&= \mathsf{ext}\big(\lambda(rf)\c \lambda(rf)\c \lambda(rrf)\big) = r(\lambda(rf))
\end{align*}
so that $L(d, \lambda f, \lambda f, \lambda(rf)) = \mathsf{subst}(r(\lambda(rf)), d(f)) =
d(f)$ as required.

It remains only to prove the Lemma. We will derive \eqref{tagged} by $\Id$-elimination on $p$,
for which it suffices to consider the case where $u = v$ and $p = r(u)$. So we must show that
\begin{equation*}
\inferrule{
  u : \Pi(A, B)\\
  c : C(\lambda x.\, u \cdot x\c \lambda x.\, u \cdot x\c \lambda x.\, r(u \cdot x))
  }
  {\phi(r(u), c) : C(u, u, \lambda x.\, r(u) \ast x)}
\end{equation*}
is derivable; which in turn we may do by $\Pi$-elimination on $u$. Indeed, when we have $u =
\lambda f$ for some $f : (x : A)\, B(x)$, we find that $C(\lambda x.\, u \cdot x\c \lambda
x.\, u \cdot x\c \lambda x.\, r(u \cdot x)) = C(\lambda f, \lambda f, \lambda(rf)) = C(u, u,
\lambda x.\, r(u) \ast x)$ so that we may take $\phi(r(\lambda f), c) = c$.
\end{proof}

\pgph We end the paper with an informal discussion of the adequacy of our
strengthening of the principle of function extensionality. We have portrayed it
as a necessary strengthening, but we have not indicated why we think it
sufficient: could there not be yet more exotic propositional equalities of the
sort considered in Examples \ref{section56} which our $\Pi$-$\Id$-elimination
rule cannot verify? The reason the author believes this not to be the case is
essentially semantic. As mentioned in~\S\ref{funextdef}, if we wish to describe
higher-dimensional categorical semantics for type theory in which $\Pi$-type
formation is a (suitably weak) right adjoint to reindexing, then we need a form
of function extensionality. As it turns out, the $\Pi$-$\Id$-elimination rule
given above is just what is needed to make this go through. The author has
verified the details of this for two-dimensional models in~\cite{g1}, and has
sketched them for a putative theory of three-dimensional models. Moreover,
there is a general argument which suggests that this extends to all higher
dimensions, which runs as follows.

When we form higher-dimensional models of type theory, we obtain the
higher dimensionality from the identity type structure. In order for $\Pi$-type
formation to provide a weak right adjoint to pullback, it must respect the
higher-dimensionality, and hence the identity type structure. Now, if we are
given $A \ty$ and $x : A \t B(x) \ty$, then dependent product formation over $x
: A$ sends the identity type
\begin{equation*}
    x : A\c y, z : B(x) \t \Id_{B(x)}(y, z) \ty
\end{equation*}
to the type
\begin{equation*}
    m, n : \Pi(A, B) \t \Pi x : A.\, \Id_{B(x)}(m \cdot x, n \cdot x)  \ty\text;
\end{equation*}
and to say that function space formation preserves the identity type structure is to say that
this latter type should act like an identity type for $\Pi(A, B)$; and it precisely this which
is expressed by our elimination rule $\Pi$-$\Id$-\textsc{elim}.

\bibliographystyle{acm}
\bibliography{biblio}

\end{document}